\documentclass[number,citesort,dvips]{arxbj}

\aid{0}
\volume{17}
\issue{1}
\pubyear{2011}
\firstpage{276}
\lastpage{289}
\doi{10.3150/10-BEJ262}

\makeatletter
\newcommand{\NN}{\mathbb N}
\newcommand{\RR}{\mathbb R}
\newcommand{\dto}{\stackrel{d}{\longrightarrow}}
\newcommand{\E}{\mathrm{E}}
\newcommand{\p}{\mathrm{P}}
\newcommand{\eps}{\varepsilon}
\newcommand{\sign}{\operatorname{sign}}
\newtheorem{theo}{Theorem}[section]
\newtheorem{prop}{Proposition}[section]
\newtheorem{lemma}{Lemma}[section]
\newtheorem{cor}{Corollary}[section]
\newproclaim{defi}{Definition}[section]
\newremark{remark}{Remark}
\newcommand{\eqref}[1]{(\ref{#1})}
\DeclareMathAlphabet\mathcaligr{OMS}{cmsy}{m}{n}
\newcommand{\cal}{\mathcaligr}
\makeatother

\begin{document}
\begin{frontmatter}

\title{On the heavy-tailedness of Student's $t$-statistic}
\runtitle{On the heavy-tailedness of Student's $t$-statistic}

\begin{aug}
\author{\fnms{Fredrik} \snm{Jonsson}\corref{}\thanksref{}\ead[label=e1]{jonsson@math.uu.se}} %\and
\address{Department of Mathematics,
Uppsala University, Box 480, SE-751 06 Uppsala, Sweden.\\
\printead{e1}}
\end{aug}

% HISTORY:
\received{\smonth{8} \syear{2009}}
\revised{\smonth{1} \syear{2010}}

% ABSTRACT
%
\begin{abstract}
Let $\{X_i\}_{i\geq1}$ be an i.i.d. sequence of random variables and
define, for $n\geq2$,
\[
T_n=
\cases{
n^{-1/2}\hat\sigma_n^{-1}S_n, &\quad$\hat\sigma_n >0$,\cr
0, &\quad$\hat\sigma_n =0$,
} \hspace*{18pt}
\mbox{with }
S_n=\sum_{i=1}^nX_i,\
\hat\sigma^2_n=\frac{1}{n-1}\sum_{i=1}^n (X_i-n^{-1}S_n )^2.
\]
We investigate the connection between the distribution of an
observation $X_i$ and finiteness of $\mathrm{E}|T_n|^r$ for $(n,r)\in
\mathbb N_{\geq2} \times\mathbb R^+$.
Moreover, assuming $T_n\dto T$, we prove that for any $r>0$, $\lim
_{n\to\infty}\mathrm{E}|T_n|^r=\mathrm{E}
|T|^r<\infty$, provided there is an integer $n_0$ such that $\mathrm
{E}|T_{n_0}|^r$ is finite.
\end{abstract}

% KEYWORDS
%

\end{frontmatter}

%s1 ###
\section{Introduction}\label{intro}
%%%%%%%%%%%%%%%%%%%%%%%%%%%%%%%%%%%%%%%%%%%%%%%%%%%%%
Assume, in the following, that $\{X_i\}_{i \geq1} $ is a sequence of
independent random variables, each with distribution $F$. Then, for
$n\geq2$, define the \textit{$t$-statistic} random variables
\[
T_n=
\cases{
n^{-1/2}\hat\sigma_n^{-1}S_n, &\quad $\hat\sigma_n >0$,\cr
0, &\quad $\hat\sigma_n =0$,
}
\hspace*{18pt}\mbox{with }
S_n=\sum_{i=1}^nX_i,\
\hat\sigma^2_n=\frac{1}{n-1}\sum_{i=1}^n (X_i-n^{-1}S_n )^2.
\]
In the case where $F$ is a normal distribution with mean zero, the
distribution of $T_n$ is the well-known \textit{$t$-distribution with
$n-1$ degrees of freedom}.
The effect of non-normality of $F$ on the distribution of $T_n$ has
received considerable attention in the statistical literature. For a
review, see \cite{Zabell}.
$t$-distributions do not only occur in the inference of means, but also
sometimes in models of data in the economic sciences; see \cite{Praetz}.
There seem to be two characteristic properties which, in comparison
with the normal distribution, make these distributions convenient in
certain modeling situations: a higher degree of heavy-tailedness
(moments are finite only below the degree of freedom) and a higher
degree of so-called kurtosis.

This paper investigates the tail behaviour of $T_n$ and the related
issue of the existence of moments $\E|T_n|^r$, for a parameter $r>0$,
under more general conditions than the normal assumption.
Motivating questions were the following:
Is it generally true that $\E|T_n|^r$ can only be finite for $r<n-1$?
For which kinds of distributions is the converse implication false?
Assuming the\vspace*{-2pt} often encountered $T_n\dto T$,
is it then generally true that $\E|T_n|^r \to\E|T|^r$?

%%%%%%%%%%%%%%%%%%%%%%%%%%%%%%%%%%%%%%
%s2 ###
\section{Summary}
The fundamental result is Theorem \ref{prop3}, which presents two
conditions, each equivalent to finiteness of $\E|T_n|^r$. The result
is based on a connection between the tail behaviour of $T_n$ and
probabilities of having almost identical observations $X_1, \ldots,X_n$.
Theorem \ref{thm-next} states that finiteness of $\E|T_n|^r$ implies
finiteness of $\E|T_{n+1}|^r$,
and is followed by Theorem \ref{thm0} which states that $t$-statistic
random variables never possess moments above the degree of freedom
unless $F$ is discrete.
It is established in Section \ref{tre}, under the assumption that $F$
is continuous, that \textit{regularity},
referring to the degree of heavy-tailedness of $t$-statistic random
variables, is measurable in terms of the behaviour of certain \textit
{concentration functions} related to $F$.
Theorem \ref{thm1} states that $\lim_{n\to\infty}\E|T_n|^r=\E
|T|^r$ whenever there is an integer $n_0$ such that $\E|T_{n_0}|^r$ is
finite and $\{T_n\}$ converges in distribution.

\begin{remark*}
This paper is an abridged version of \cite{jonsson}.
The results found in Section \ref{tre} here are there generalized
beyond the continuity assumption.
We also refer to \cite{jonsson} for a discussion of related results
previously obtained by H.~Hotelling.
\end{remark*}
%
%%%%%%%%%%%%%%%%%%%%%%%%%%%%%%%%%%%%%%%%%%%%%%%%%%%%%

%s3 ###
\section{Characterizing $\E|T_n|^r<\infty$ through bounds on $\p
(|T_n|>x)$} \label{ett}
%%%%%%%%%%%%%%%%%%%%%%%%%%%%%%%%%%%%%%%%%%%%%%%%%%%%%
A close connection exists between $T_n$ and the \textit
{self-normalized sum} $S_n/V_n$; see Lemma \ref{lem00} (whose
elementary proof we omit).
The connection allows $\E|T_n|^r$ to be expressed with probabilities
relating to $S_n/V_n$, as in Lemma \ref{lem0},
revealing that finiteness of $\E|T_n|^r$ depends on the magnitude of
the probabilities of having $S_n/V_n$ close to $\pm\sqrt{n}$. Some
geometric relations between $S_n/V_n$ close to $\pm\sqrt{n}$ and
almost identical observations $X_1, \ldots,X_n$ are then given in
Lemmas \ref{lem1} and~\ref{lem2}.

%%%%%%%%%%%%%%%%%%%%%%%%%%%%%%%%%%%%%%%%%%%%%%%%%%%%%%%%%%%%%%%%
%
\begin{lemma} \label{lem00}
Define
\[
V_n= \Biggl(\sum_{i=1}^nX_i^2 \Biggr)^{1/2},\qquad U_n^*=
\cases{
0, &\quad $S_n/V_n=n$ or $V_n=0$,\cr
(S_n/V_n)^2, &\quad otherwise.
}
\]
It then holds, for any $x\geq0$, that $T_n^2> x$ if and only if
$U_n^*>nx/(n+x-1)$.
\end{lemma}

%%%%%%%%%%%%%%%%%%%%%%%%%%%%%%%%%%%%%%%%%%%%%%5
%
\begin{lemma} \label{lem0}
For $r>0$ and $U_n^*$ as in Lemma \ref{lem00},
\[
\E|T_n|^r=\frac{r}{2}n(n-1)^{r/2}\int_{0}^n z^{r/2-1}\p
(U_n^*>z )
(n-z)^{-(r/2+1)}\,\mathrm{d}z.
\]
\end{lemma}
%
%%%%%%%%%%%%%%%%%%%%%%%%%%%%%%%%%%%%%%%%%%%%%%5

%%%%%%%%%%%%%%%%%%%%%%%%%%%%%%%%%%%%%%%%%%%%%%%%%%%%%
%
\begin{lemma} \label{lem1}
Let $\mathbf{x}=(x_1,\ldots,x_n)\in\RR^n$ and $h\in(0,1)$ be given
such that $x_1\neq0$ and $n-u_n < h^2$ with
$u_n=(\sum_{i=1}^nx_i)^2/\sum_{i=1}^nx_i^2$. Then, with $C_1=\sqrt{5}$,
\[
|x_i-x_1| < h C_1|x_1|\qquad \mbox{for all }i\neq1.
\]
Moreover, $C_1=C_1(n,h)=\sqrt{2+2h+h^2}$
is optimal for the conclusion to be valid for all $\mathbf{x}$.
\end{lemma}

%%%%%%%%%%%%%%%%%%%%%%%%%%%%%%%%%
%
\begin{lemma} \label{lem2}
Let $\mathbf{x}=(x_1,\ldots,x_n)\in\RR^n$ and $h\in(0,1)$ be given
such that, with $C_2=1$,
\[
|x_i-x_1| < C_2h|x_1|/\sqrt{n-1} \qquad\mbox{for all }i\neq1.
\]
Then $n-u_n < h^2$ with
$u_n=(\sum_{i=1}^nx_i)^2/\sum_{i=1}^nx_i^2$.
Moreover, in the case where $n$ is odd, $C_2=C_2(n,h)$ must satisfy
$C_2\leq\sqrt{n/(n-h^2)}$
for the conclusion to be valid for all $\mathbf{x}$.
\end{lemma}
%
%%%%%%%%%%%%%%%%%%%%%%%%%%%%%%%%%%%%%%%%%%%%%%%%%%%%%%%%%%%%
%
\begin{theo} \label{prop3}
The following three quantities are either all finite or all infinite:
\begin{eqnarray*}
\mbox{\textup{(i)}}&&\hspace*{5pt} \E|T_n|^r ; \\
\mbox{\textup{(ii)}}&&\hspace*{5pt} \E\Biggl(|X_1|^r\bigwedge_{i=2}^n|X_i-X_1|^{-r} I \{
|X_i-X_1|>0, \mbox{ some } i\leq n \} \Biggr) ;\\
\mbox{\textup{(iii)}}&&\hspace*{5pt}\int_{x\neq0}\int_{0}^1 h^{-(r+1)} \bigl(\bigl (\p(|X-x|<
h|x| ) \bigr)^{n-1}-p_x^{n-1} \bigr)\,\mathrm{d}h \, \mathrm{d}F(x)
\qquad\mbox{with } p_x=\p(X=x) .
\end{eqnarray*}
\end{theo}

%%%%%%%%%%%%%%%%%%%%%%%%%%%%%%%%%%%%%%%%%%%%%%%%%%%%%%%%%%%%
%%%%%%%%%%%%%%%%%%%%%%%%%%%%%%%%%%%%%%%%%%%%%%%%%%%%%%%%%%%%
%
\begin{pf*}{Proof of Lemma \ref{lem0}}
By \cite{Gut}, Theorem 12.1, Chapter 2, together with Lemma \ref
{lem00} and a change of variables,
we have
\begin{eqnarray*}
\E|T_n|^r&=&\frac{r}{2}\int_0^\infty y^{r/2-1}\p(T_n^2>y)\,\mathrm
{d}y\\
&=&\frac{r}{2}\int_0^\infty y^{r/2-1}\p\bigl(U_n^*>ny/(n+y-1)
\bigr)\,\mathrm{d}y\\
&=&\frac{r}{2}n(n-1)^{r/2}\int_{0}^n z^{r/2-1}\p(U_n^*>z )
(n-z)^{-(r/2+1)}\,\mathrm{d}z.
\end{eqnarray*}
\upqed
\end{pf*}

%%%%%%%%%%%%%%%%%%%%%%%%%%%%%%%%%%%%%%%%%%%%%%%%%%%%%%%%%%%%
%
\begin{pf*}{Proof of Lemma \ref{lem1}}
We argue by contraposition.
Due to the invariance with respect to scaling of $\mathbf{x}$ and
permutation of the coordinates $x_2,\ldots,x_n$, it suffices to prove that
\[
% \label{obj}
|x_2-x_1|\geq h|x_1|\quad \Longrightarrow\quad  n-u_n \geq h^2/C_1^2
\]
with $C_1=\sqrt{2+2h+h^2}$ and that equalities are simultaneously attained.
Set $x_2=x_1+\eps$ and $\underline{x} =(x_3, \ldots, x_n)$. We then
minimize $n-u_n$ with respect to $\underline{x}$ and $\eps$.
Note that
%
%e1 ###
\begin{equation} \label{fderiv}
\frac{\partial(n-u_n)}{\partial x_j} = \frac{-2\sum_{i=1}^n x_i
(\sum_{i=1}^n x_i^2-x_j\sum_{i=1}^n x_i )}{ (\sum_{i=1}^n
x_i^2 )^2}.
\end{equation}
First, set \eqref{fderiv} to zero for $j=3, \ldots,n$. Since
$\sum x_i=0$ corresponds to $u_n=0$, which is non-interesting with
respect to the minimization of $n-u_n$,
these equations reduce to
%
%e2 ###
\begin{equation} \label{lem1a}
\sum_{i=3}^n x_i^2-x_j \sum_{i=3}^nx_i = x_j(x_1+x_2) - (x_1^2+x_2^2)
\qquad\mbox{for } j=3,\ldots,n.
\end{equation}
We claim that (\ref{lem1a}) has the unique solution
%
%e3 ###
\begin{equation} \label{lem1b}
x_j = (x_1^2+x_2^2)/(x_1+x_2) = (2x_1^2+2x_1\eps+\eps^2)/(2x_1+\eps
)\qquad\mbox{for } j=3,\ldots,n.
\end{equation}
To verify this, assume that $\underline x$ is a solution of (\ref{lem1a}).
Since $\sum_{i=3}^n x_i^2$ and $\sum_{i=3}^n x_i$ do not vary with
$j$, $\underline x$ must be of the form $x_j=\mathit{const}.$, $j=3,\ldots,n$.
However, the left-hand side of (\ref{lem1a}) then vanishes for all
$j$, which gives (\ref{lem1b}) as the unique solution.
Inserting the solution into $n-u_n$ gives
%
%e4 ###
\begin{equation} \label{lem1c}
(n-u_n)_{\min}(\eps) = \eps^2/(x_1^2+x_2^2) = \eps
^2/(2x_1^2+2x_1\eps+\eps^2).
\end{equation}
It remains to minimize with respect to $\eps$ with $\eps\notin
(-h|x_1|,h|x_1|)$.
The equation
\[
\frac{\partial}{\partial\eps} \biggl(\frac{\eps^2}{2x_1^2+2x_1\eps
+\eps^2} \biggr)=0
\]
has the unique solution $\eps=-2x_1$ which cannot be a minimum since
a minimum must satisfy $\sign(\eps)=\sign(x_1)$, by the
representation (\ref{lem1c}).
The solution is hence obtained for $\eps=\sign(x_1) h|x_1|$,
\[
(n-u_n)_{\min}=(hx_1)^2/\bigl(x_1^2(2+2h+h^2)\bigr)=h^2/(2+2h+h^2).
\]
It follows that $C_1=C_1(h)=\sqrt{2+2h+h^2}\leq\sqrt{5}$ is an
optimal constant, as claimed.
\end{pf*}

%%%%%%%%%%%%%%%%%%%%%%%%%%%%%%%%%%%%%%%%%%%%%%%%%%%%%%%%%%%%
%
\begin{pf*}{Proof of Lemma \ref{lem2}}
Assume that
%
%e5 ###
\begin{equation} \label{lem2'''a}
|x_i-x_1| < C_2h|x_1|/\sqrt{n-1} \qquad\mbox{for all } i=2,\ldots, n.
\end{equation}
The aim is to verify that $n-u_n<h^2$
with $C_2=C_2(n,h)$ optimally large. We therefore maximize $n-u_n$
over the rectangular region \eqref{lem2'''a} with $x_1\neq0$, $C_2$
and $h$ fixed.
It suffices to consider the restriction of $n-u_n$ to the corners of
the region \eqref{lem2'''a} since the maximum attained at a point
$y=(y_1,\ldots,y_n)$ in the interior of the region, or in the interior
of an edge,
would mean that, for some $j=2,\ldots, n$ and some $\eta>0$,
%
%e8 ###
%e7 ###
%e6 ###
\begin{eqnarray}\label{lem2b}
\frac{\partial(n-u_n)}{\partial x_j}(y)&=&0,\\\label{lem2c}
\frac{\partial(n-u_n)}{\partial x_j} (y_1,\ldots
,y_{j-1},y_j-h,y_{j+1},\ldots, y_n )&\geq&0 \qquad\mbox{for all }
0<h<\eta,\\\label{lem2d}
\frac{\partial(n-u_n)}{\partial x_j} (y_1,\ldots
,y_{j-1},y_j+h,y_{j+1},\ldots, y_n )&\leq&0 \qquad\mbox{for all }
0<h<\eta.
\end{eqnarray}
Recall, from the proof of Lemma \ref{lem1}, that
\[
\frac{\partial(n-u_n)}{\partial x_j} = \frac{-2\sum_{i=1}^n x_i
(\sum_{i\neq j} x_i^2-x_j\sum_{i\neq j} x_i )}{ (\sum
_{i=1}^n x_i^2 )^2}.
\]
We may assume that $C_2h<\sqrt{n-1}$
since the point $x_i\equiv0$ would otherwise belong to the region
yielding $u_n=1$, in which case $n-u_n<h^2$ cannot hold.
This implies that $\sign(x_i)=\sign(x_1)$ for all $i=2,\ldots, n$
so that neither $\sum x_i$ nor $\sum_{i\neq j} x_i$ change sign within
the region.
Assume, due to invariance with respect to scaling, that $x_1>0$.
Conditions \eqref{lem2b}--\eqref{lem2d} may then be reformulated as
\[
\sum_{i\neq j} y_i^2-y_j\sum_{i\neq j} y_i=0,\qquad
\sum_{i\neq j} y_i^2-(y_j-h)\sum_{i\neq j} y_i<0,\qquad
\sum_{i\neq j} y_i^2-(y_j+h)\sum_{i\neq j} y_i>0,
\]
which is contradictory since $h>0$ and $\sum_{i\neq j} y_i>0$.

Now, consider the restriction of $n-u_n$ to the corners of the region
\eqref{lem2'''a}.
Set $k:=| \{i\dvtx x_i=x_1+\eps\}|-| \{i\dvtx x_i=x_1-\eps\}
|$ so that
%
%e9 ###
\begin{eqnarray}\label{lem2a}
n-u_n
&=&\frac{n(nx_1^2+(n-1)\eps^2+2k\eps x_1)-(nx_1+k\eps
)^2}{nx_1^2+(n-1)\eps^2+2k\eps x_1} \nonumber
\\[-8pt]
\\[-8pt]
&=&\frac{\eps^2(n(n-1)-k^2)}{nx_1^2+(n-1)\eps^2+2k\eps x_1} = \frac
{h^2C_2^2(n-k^2/(n-1))}{n+C_2^2h^2+2kC_2h/\sqrt{n-1}}.\nonumber
\end{eqnarray}
Take $C_2=1$ in (\ref{lem2a}) and $z=k(n-1)^{-1/2}$. Algebraic
manipulations yield
\[
\frac{n-k^2/(n-1)}{n+h^2+2kh/\sqrt{n-1}}\leq1 \quad \Longleftrightarrow\quad
(h+z)^2\geq0
\]
so that $C_2=1$ is sufficiently small for the desired bound $n-u_n<h^2$.
We find, by taking $k=0$ in (\ref{lem2a})  (which is possible when $n$
is odd) that
\[
C_2^2n/(n+C_2^2h^2)\leq1 \quad \Longleftrightarrow \quad C_2^2\leq n/(n-h^2)
\]
so that $C_2\leq\sqrt{n/(n-h^2)}$
is then necessary for $n-u_n<h^2$ to hold.
\end{pf*}
%
%%%%%%%%%%%%%%%%%%%%%%%%%%%%%%%%%%%5

%%%%%%%%%%%%%%%%%%%%%%%%%%%%%%%%%%%%%%%%%%%%%%%%%%%%%
%%%%%%%%%%%%%%%%%%%%%%%%%%%%%%%%%%%%%%%%%%%%%%%%%%%%%%%%%%%%
%
\begin{pf*}{Proof of Theorem \ref{prop3}}
We first deduce the equivalence between (i) and (iii).
By Lemma \ref{lem0}, we find that $E|T_n|^r<\infty$
is equivalent to, for some $\delta<1$,
\[
\int_{n-\delta}^n z^{r/2-1}\p(U_n^*>z )
(n-z)^{-(r/2+1)}\,\mathrm{d}z<\infty\quad \Longleftrightarrow\quad \int
_{0}^\delta h^{-(r+1)}\p(n-U_n^*<h^2 )\,\mathrm{d}h<\infty,
\]
which, in turn, is equivalent to
%
%e10 ###
\begin{equation}\label{prop3_1}
\int\!\!\int_{0}^\delta h^{-(r+1)}\p(0<n-U_n<h^2\mid X_1=x )\,
\mathrm{d}h\,\mathrm{d}F(x)<\infty.
\end{equation}
The event $X_1=0$ implies $U_n\leq n-1$ by the Cauchy--Schwarz
inequality so that \eqref{prop3_1} reduces to
\[
\int_{x\neq0} \int_{0}^\delta h^{-(r+1)}\p(0<n-U_n<h^2\mid
X_1=x )\,\mathrm{d}h\,\mathrm{d}F(x)<\infty,
\]
which is equivalent to
\[
\int_{x\neq0} \int_{0}^\delta h^{-(r+1)}\p(n-U_n<h^2\mid
X_1=x )-p_x^{n-1}\,\mathrm{d}h\,\mathrm{d}F(x)<\infty
\]
since $U_n=n$ corresponds to $X_i=X_1$ with $p_x=\p(X=x)$.
Finally, apply Lemmas \ref{lem1} and \ref{lem2}, and set $\delta=1$
to arrive at condition (iii).

For the equivalence between (ii) and (iii), define
$ A_n= \{|X_i-X_1|>0$, some $i \leq n \}$.
Condition on $X_1$ and convert expectation into integration of tail
probabilities (cf. \cite{Gut}, Theorem 12.1, Chapter~2): %Lemma
\begin{eqnarray*}
\E\Biggl(|X_1|^r\bigwedge_{i=2}^n|X_i-X_1|^{-r} I_{A_n} \Biggr)
&=&\int_{x \neq0} \E\Biggl(\bigwedge_{i=2}^n
(|X_i-x||x|^{-1} )^{-r} I_{A_n} \Biggr)\,\mathrm{d}F(x)\\
&=&
r\int_{x \neq0} \int_0^\infty h^{-(r+1)}  \bigl(\p(|X-x|< h
|x| \bigr)^{n-1}-p_x^{n-1}\,\mathrm{d}h\,\mathrm{d}F(x).
\end{eqnarray*}
The equivalence between (ii) and (iii) then follows from the fact that
\begin{eqnarray*}
&&\int_{x \neq0} \int_1^\infty h^{-(r+1)} \bigl(\p(|X-x|< h
|x| ) \bigr)^{n-1}\,\mathrm{d}h \, \mathrm{d}F(x)\\
&&\quad \leq\int_{x \neq0} \int_1^\infty h^{-(r+1)}\,\mathrm{d}h \,
\mathrm{d}F(x) < \infty.
\end{eqnarray*}
\upqed
\end{pf*}

%%%%%%%%%%%%%%%%%%%%%%%%%%%%%%%%%%%%%%%%%%%%%%%%%%%%%
%s4 ###
\section{Two general facts regarding finiteness of $\E
|T_{n}|^r$}\label{tva}
%%%%%%%%%%%%%%%%%%%%%%%%%%%%%%%%%%%%%%%%%%%%%%%%%%%%%

\begin{theo} \label{thm-next}
For any couple $(n,r)\in\NN_{\geq2} \times\RR^+$, if $\E
|T_{n}|^r$ is finite, then so is
$\E|T_{n+1}|^r$.
\end{theo}
%
%%%%%%%%%%%%%%%%%%%%%%%%%%%%%%
%
\begin{pf}
Due to Theorem \ref{prop3}, it suffices to show that
%
%e11 ###
\begin{equation} \label{6a}
\E\Biggl[|X_1|^r\bigwedge_{i=2}^n|X_i-X_1|^{-r} I_{A_n} \Biggr]<\infty
\quad \Longrightarrow\quad
\E\Biggl[|X_1|^r\bigwedge_{i=2}^{n+1}|X_i-X_1|^{-r} I_{A_{n+1}}
\Biggr]<\infty,
\end{equation}
where $ A_k:= \{|X_i-X_1|>0$, some $i \leq k \}$.
Define $ A'_n= \{|X_i-X_1|>0$, some $3\leq i \leq n+1 \}$.
It follows that $A_{n+1} = A_n \cup A'_n$
so that
$ I_{A_{n+1}} \leq I_{A_{n}} + I_{A'_{n}}$, which gives
\begin{eqnarray*}
&&\E\Biggl[|X_1|^r\bigwedge_{i=2}^{n+1}|X_i-X_1|^{-r} I_{A_{n+1}} \Biggr]
\\
&&\quad  \leq\E\Biggl[|X_1|^r\bigwedge_{i=2}^{n+1}|X_i-X_1|^{-r}
I_{A_{n}} \Biggr] +
\E\Biggl[|X_1|^r\bigwedge_{i=2}^{n+1}|X_i-X_1|^{-r} I_{A'_{n}} \Biggr]
\\
&&\quad  \leq\E\Biggl[|X_1|^r\bigwedge_{i=2}^{n}|X_i-X_1|^{-r}
I_{A_{n}} \Biggr] +
\E\Biggl[|X_1|^r\bigwedge_{i=3}^{n+1}|X_i-X_1|^{-r} I_{A'_{n}} \Biggr]\\
&&\quad  = 2\E\Biggl[|X_1|^r\bigwedge_{i=2}^{n}|X_i-X_1|^{-r} I_{A_{n}} \Biggr].
\end{eqnarray*}
The conclusion follows.
\end{pf}

%%%%%%%%%%%%%%%%%%%%%%%%%%%%%%%%%%%%%%%%%%%%%%%%%
%
\begin{theo} \label{thm0}
Assume that $F$ decomposes into $F_d+F_c$, with discrete and continuous
measures~$F_d$ and $F_c$, respectively, and that $F_c \not\equiv0$.
It is then necessary that $r<n-1$ for $\E|T_n|^r$ to be finite.
\end{theo}
\begin{pf}
Let $F_c$ have total mass $\eps>0$. It suffices to verify that $\E
|T_n|^{n-1}$ is infinite,
which, by Theorem~\ref{prop3}, is equivalent to
\[
\int_{x \neq0} \int_0 ^1 h^{-n}\bigl ( \bigl(\p(|X-x|< h |x|
) \bigr)^{n-1}-p_x^{n-1} \bigr)\,\mathrm{d}h\,\mathrm{d}F(x)=\infty.
\]
The last identity is a consequence of
%
%e12 ###
\begin{equation}\label{p3a}
\int\!\!\int_0 ^1 h^{-n} \bigl(\p(|X-x|< h |x| )
\bigr)^{n-1}\,\mathrm{d}h\,\mathrm{d}F_c(x)=\infty.
\end{equation}
To verify \eqref{p3a}, consider the restriction of $F_c$ to a set
$[-C,-1/C] \cup[1/C,C]$ with $C$ sufficiently large so that the
restricted measure still has positive mass.
It then suffices to establish the condition
%
%e13 ###
\begin{equation}\label{p3d}
\hspace*{-5pt}\int \bigl(\p(|X-x|< h )h^{-1} \bigr)^{n-1}\,\mathrm
{d}F_c(x)>\eta_n \qquad\mbox{for all $h$ and some constant } \eta_n=\eta
_n(F_c,n).
\end{equation}
First, consider $n=2$. Discretize $[-C,C]$ uniformly with interval
length $h$, that is, put
$x_k=hk$ for $k\in[-N,N]$ and $N= \lceil Ch^{-1}\rceil$.
Then
\begin{eqnarray*}
\int\p(|X_c-x|< h )\,\mathrm{d}F_c(x)&=&\sum
_{k=-N}^{k=N}\int_{x_{k-1}}^{x_k}\p(|X_c-x|< h )\,\mathrm
{d}F_c(x)  \\
&\geq& \sum_{k=-N}^{k=N}\int_{x_{k-1}}^{x_k}\p\bigl(X_c\in
(x_{k-1},x_k] \bigr)\,\mathrm{d}F_c(x)  \\
&=&\sum_{k=-N}^{k=N} \bigl(\p\bigl(X_c\in(x_{k-1},x_k] \bigr) \bigr)^2.
\end{eqnarray*}
Applying the Cauchy--Schwarz inequality, we obtain
\[
\sum_{k=-N}^{k=N} \bigl(\p(X_c\in(x_{k-1},x_k] ) \bigr)^2\geq
\Biggl(\sum_{k=-N}^{k=N}\p\bigl(X_c\in(x_{k-1},x_k] \bigr) \Biggr)^2 (2N)^{-1}
= \eps^2 (2N)^{-1} \geq C^{-1}\eps^2h.
\]
Conclusion (\ref{p3d}) follows with $\eta_2=C^{-1}\eps^2$.
For $n>2$, an application of the H\"older inequality yields
\[
\eta_2^{n-1} \leq \biggl(\int\p(|X_c-x|< h)h^{-1}\,\mathrm
{d}F_c(x) \biggr)^{n-1}
\leq\eps^{n-2}\int\bigl(\p(|X_c-x|< h)h^{-1} \bigr)^{n-1}\,\mathrm{d}F_c(x).
\]
The desired conclusion (\ref{p3d})
follows with $\eta_n=\eta_2^{n-1}\eps^{2-n}$.
\end{pf}

%%%%%%%%%%%%%%%%%%%%%%%%%%%%%%%%%%%%%%%%%%%%%%%%%%%%%

%s5 ###
\section{Regularity and concentration functions}\label{tre}
%%%%%%%%%%%%%%%%%%%%%%%%%%%%%%%%%%%%%%%%%%%%%%%%%%%%%
%
\begin{defi} \label{def2}
Given the distribution of a random variable $X$, define the
concentration functions~$q$ and $Q$,
for real-valued arguments $h\geq0$, by
\[
Q(h) = \sup_x {\p(|X-x| \leq h )},\qquad
q(h) = \sup_x {\p(|X-x| \leq|x|h )}.
\]
\end{defi}

$Q$ is known as the \textit{L\'evy concentration function}.
Theorem \ref{thm3} below characterizes finiteness of $\E|T_n|^r$ in
terms of the limiting behaviour of $q(h)$ as $h$ tends to zero.
Note that a statement of the kind ``$Q(h)={\cal O}(h^{\lambda})$''
(for some $\lambda\leq1$) refers to the \textit{local} behaviour of
the distribution.
The most regular behaviour in this respect is that of an absolutely
continuous distribution with bounded density function, in which case
$ Q(h)={\cal O}(h)$, while $\lambda<1$ typically corresponds to one or
several ``explosions'' of the density function. The Cantor
distributions also form fundamental examples of such irregularity (cf.
\cite{jonsson}, pages 29--31).
The parameter $\lambda$ has, in this sense, a meaning of ``degree of
irregularity'' concerning the distribution, with smaller values of
$\lambda$ indicating higher degrees of irregularity. A statement
$q(h)={\cal O}(h^{\lambda})$, on the other hand, also has a \textit
{global} component. It requires more regularity of the distribution
``at infinity'' compared with $Q(h)={\cal O}(h^{\lambda})$, while, at
the same time, being
less restrictive regarding the local behaviour of the distribution at
the origin.

%%%%%%%%%%%%%%%%%%%%%%%%%%%%%%%%%%%%%%%%%%%%%%%%%%%%%
%
\begin{theo} \label{thm3}
The following two implications hold for any continuous probability
measure $F$:
\begin{eqnarray*}
\mbox{\textup{(i)}}&&\hspace*{5pt} q(h)={\cal O}(h^{\lambda}) \mbox{ for some } \lambda>r/(n-1)
\Longrightarrow\E|T_n|^r<\infty;\\
\mbox{\textup{(ii)}}&&\hspace*{5pt} \E|T_n|^r<\infty\Longrightarrow q(h)={\cal O}(h^{\lambda
})\mbox{ with } \lambda=r/n .
\end{eqnarray*}
\end{theo}\vspace*{-1.5pt}
%
%%%%%%%%%%%%%%%%%%%%%%%%%%%%%%%%%%%%%%%%%%%%%%%%%%%%%%%%%

A simple criterion guaranteeing the optimal $q(h)={\cal O}(h)$ is given
by the following proposition.\vspace*{-1.5pt}

%%%%%%%%%%%%%%%%%%%%%%%%%%%%%%%%%%%%%%%%%%%%%%%%%%%%%%%%%
%
\begin{prop} \label{prop4}
The property $q(h)={\cal O}(h)$ is obtained for any absolutely
continuous distribution $F$ with bounded density function $f$
satisfying the assumption of a positive constant $N$ such that
%
%e14 ###
\begin{equation}\label{prop4villk}
f(x_2)\leq f(x_1) \qquad\mbox{for any } x_1, x_2 \mbox{ such that }
N\leq x_1\leq x_2 \mbox{ or } -N\geq x_1\geq x_2.
\end{equation}
\end{prop}
%
%%%%%%%%%%%%%%%%%%%%%%%%%%%%%%%%%%%%%%%%%%%%%%%%%%%%

%%%%%%%%%%%%%%%%%%%%%%%%%%%%%%%%%%%%%%%%%%%%%%%%%%%%
%
\begin{pf*}{Proof of Theorem \ref{thm3}}
For (i), condition (iii) of Theorem \ref{prop3} reads, by continuity,
%
%e15 ###
\begin{equation}\label{thm3''}
\int_{x\neq0}\int_{0}^1 h^{-(r+1)} \bigl(\p(|X-x|< h|x|
) \bigr)^{n-1}\,\mathrm{d}h \, \mathrm{d}F(x)<\infty.
\end{equation}
Applying the assumption on $q$ to the integrand yields
\begin{eqnarray*}
&&\int_{x\neq0}\int_{0}^1 h^{-(r+1)} \bigl(\p(|X-x|< h|x|
) \bigr)^{n-1}\,\mathrm{d}h \, \mathrm{d}F(x)\\
&&\quad \leq C\int_{x\neq0}\int_{0}^1 h^{-(r+1)}h^{\lambda(n-1)}\,\mathrm
{d}h \, \mathrm{d}F(x) = C\int_{0}^1 h^{-(r+1)}h^{\lambda
(n-1)}\,\mathrm{d}h\\
&&\quad = C/\bigl(\lambda(n-1)-r\bigr),
\end{eqnarray*}
which proves (\ref{thm3''}).
To verify the second implication, we argue by contraposition. Assume that
%
%e16 ###
\begin{equation}
q(h)\not={\cal O}(h^{\lambda}) \qquad\mbox{with } \lambda=r/n. \label{thm3'}
\end{equation}
It suffices, by condition (ii) of Theorem \ref{prop3} and the
assumption of continuity, to prove that
%
%e17 ###
\begin{equation}
\E\Biggl(|X_1|^r\bigwedge_{i=2}^n|X_i-X_1|^{-r} \Biggr)=\infty. \label{thm3a'}
\end{equation}
Statement (\ref{thm3'}) is equivalent to the existence of sequences $\{
x_k\}_{k\geq1}$ and $\{h_k\}_{k\geq1}$ such that
%
%e18 ###
\begin{equation}\label{thm3a0}
1/2>h_k>0,\qquad \lim_{k\to\infty}h_k= 0,\qquad
\lim_{k\to\infty} h_k^{-r/n}\p(|X-x_k| \leq|x_k|h_k ) =
\infty.
\end{equation}
Define intervals $I_k= (x_k-|x_k|h_k, x_k+|x_k|h_k )$.
It then follows that for some $K$ and all $k\geq K$,
\begin{eqnarray*}
\E\Biggl(|X_1|^r\bigwedge_{i=2}^n|X_i-X_1|^{-r} \Biggr)
&\geq&\E\Biggl(|X_1|^r\bigwedge_{i=2}^n|X_i-X_1|^{-r} I \{ X_i
\in I_k, \mbox{ all }i \} \Biggr) \\
&\geq&2^{-1}|x_k|^r \E\Biggl(\bigwedge_{i=2}^n|X_i-X_1|^{-r} I
\{ X_i \in I_k, \mbox{ all }i \} \Biggr) \\
&\geq&2^{-(r+1)} |x_k|^r h_k^{-r} |x_k|^{-r}\E(I \{ X_i
\in I_k, \mbox{ all }i \} ) \\
&= &2^{-(r+1)}h_k^{-r} \bigl(\p(|X-x_k| \leq|x_k|h_k ) \bigr)^n.
\end{eqnarray*}
We conclude from (\ref{thm3a0}) that (\ref{thm3a'}) holds.
\end{pf*}
%
%%%%%%%%%%%%%%%%%%%%%%%%%%%%%%%%%%%%%%%%%%%%%%%%%%%%%%%%
%%%%%%%%%%%%%%%%%%%%%%%%%%%%%%%%%%%%%%%%%%%%%%%%%%%%
%
\begin{pf*}{Proof of Proposition \ref{prop4}}
It follows that, for $x>N$,
\[
f(x)(x-N) \leq\int_N^xf(y)\,\mathrm{d}y \leq1,\qquad f(-x)(x-N) \leq
\int_{-x}^{-N}f(y)\,\mathrm{d}y \leq1,
\]
so that $f(x)|x|\leq C$. Consequently, assuming that $x>2N$ and $h\leq1/2$,
we have
%
%e19 ###
\begin{equation}\label{prop4a}
\p(|X-x| \leq|x|h ) = \int_{|x|(1-h)}^{|x|(1+h)}f(y)\,\mathrm
{d}y \leq\frac{2C}{|x|}\int_{|x|(1-h)}^{|x|(1+h)}\,\mathrm{d}y =
4Ch.
\end{equation}
Regarding $0\leq x \leq2N$, we use the fact that $f$ is bounded,
$f\leq M$, so that
%
%e20 ###
\begin{equation}
\p(|X-x| \leq|x|h ) = \int_{|x|(1-h)}^{|x|(1+h)}f(y)\,\mathrm{d}y
\leq M \int_{2N(1-h)}^{2N(1+h)}\mathrm{d}y = 4MNh.\label{prop4b}
\end{equation}
Bounds analogous to \eqref{prop4a} and \eqref{prop4b} follow for
negative $x$, which proves that $q(h)={\cal O}(h)$.
\end{pf*}

%%%%%%%%%%%%%%%%%%%%%%%%%%%%%%%%%%%%%%%%%%%%%%%%%%%%%
%s6 ###
\section{Convergence}\label{fyra}
%%%%%%%%%%%%%%%%%%%%%%%%%%%%%%%%%%%%%%%%%%%%%%%%%%%%%
Convergence in distribution of $\{T_n\}$ to a random variable $T$
(e.g., standard normally distributed) is, due to Lemma \ref{lem0},
equivalent to convergence of $\{S_n/V_n\}$ to $T$.
A complete classification in terms of possible limit distributions with
corresponding conditions on $F$ was given recently by Chistyakov and
G\" otze (see \cite{ChistG}). The following interesting property was
derived somewhat earlier by Gin\'e, G\"otze and Mason in \cite{G-G-M}.

%%%%%%%%%%%%%%%%%%%%%%%%%%%%%%%%%%%%%%%%%%%%%%%%%%%%%
%
\begin{theo} \label{thm_G-G-M}
Let a distribution $F$ be given such that $S_n/V_n \to^d T$.
The sequence $\{S_n/V_n\}$ is then \textit{sub-Gaussian}, in the sense
that, for some constant $C$,
$\sup_n \E[\exp{ (tS_n/V_n )} ]\leq2 \exp{(Ct^2)}$.
\end{theo}

%%%%%%%%%%%%%%%%%%%%%%%%%%%%%%%%%%%%%%%%%%%%%%%%%%%%%
%
\begin{cor}\label{cor_G-G-M}
For any $F$ satisfying the condition of Theorem \ref{thm_G-G-M} with
respect to a random variable $T$ and any $r>0$, $\lim_{n\to\infty}\E
|S_n/V_n|^r=\E|T|^r<\infty$.
\end{cor}

%%%%%%%%%%%%%%%%%%%%%%%%%%%%%%%%%%%%%%%%%%%%%%%%%%%%%%%%%%%%%%%%%%%
%
\begin{pf}
The result follows from Theorem \ref{thm_G-G-M} and general properties
of integration; see, for example, \cite{Gut}, Theorem 5.9, Chapter 5,
or \cite{Gut}, Corollary 4.1, Chapter 5.
\end{pf}
We are now ready for the main result of this section.
%%%%%%%%%%%%%%%%%%%%%%%%%%%
%
\begin{theo} \label{thm1}
Let $F$, $T$ and $r$ be given as in Corollary \ref{cor_G-G-M}.
If $\E|T_{n_0}|^r$ is finite
for some $n_0\geq2$, then $\lim_{n\to\infty} \E|T_n|^r = \E|T|^r$.
\end{theo}
%
%%%%%%%%%%%%%%%%%%%%%%%%%%%

%%%%%%%%%%%%%%%%%%%%%%%%%%%%%%%%%%%%%%%%%%%%%%%%%%%%%%%%%%%%%%%%%%%
%
\begin{pf}
The case ``$X=constant$'', which leads to $T_n\equiv0$, is degenerate
and is henceforth excluded.
Recall, from Lemma \ref{lem0}, that
\[
\E|T_n|^r=\frac{r}{2}n(n-1)^{r/2}\int_{0}^n z^{r/2-1}\p
(U_n^*>z )
(n-z)^{-(r/2+1)}\,\mathrm{d}z.
\]
We split the desired conclusion $\lim_{n\to\infty} \E|T_n|^r = \E
|T|^r$ into the two conditions
%
%e22 ###
%e21 ###
\begin{eqnarray}\label{6'}
\hspace*{-17pt}\lim_{n\to\infty} \frac{r}{2}n^{r/2+1}\int_{0}^{n-\delta
}z^{r/2-1} \p(U_n^*>z )
(n-z)^{-(r/2+1)}\,\mathrm{d}z &=& \E|T|^r \quad\mbox{for any } 0<\delta
<1,\\\label{6}
\hspace*{-17pt}\lim_{n\to\infty} n^{r}\int_{n-\delta}^n \p(U_n^*>z )
(n-z)^{-(r/2+1)}\,\mathrm{d}z &=& 0\quad\mbox{for some } 0<\delta<1.
\end{eqnarray}
Replace \eqref{6}, via a change of variables $n-z=h^2$, by the condition
\[
\lim_{n\to\infty} n^{r}\int_0^{\delta} h^{-(r+1)} \p
(n-U_n^*<h^2 )
\,\mathrm{d}z = 0\qquad\mbox{for some } 0<\delta<1,
\]
which, in turn, by the same steps as in the proof of Theorem \ref
{prop3}, we find to be equivalent to
%
%e23 ###
\begin{eqnarray} \label{'6}
&&\lim_{n\to\infty} R_{n,\delta}=0,\nonumber
\\[-8pt]
\\[-8pt]
&&\quad R_{n,\delta}:=\int_{x\neq
0}\int_{0}^\delta n^rh^{-(r+1)} \bigl( \bigl(\p(|X-x|< h |x|
) \bigr)^{n-1}-p_x^{n-1} \bigr)\,\mathrm{d}h \, \mathrm{d}F(x)\nonumber
\end{eqnarray}
for some $0<\delta<1$ (with $p_x=\p(X=x)$). We separate the
verifications of \eqref{6'} and \eqref{'6} into Lemmas \ref{lem4} and
  \ref{lem3}, respectively.
Note that the assumption $\E|T_{n_0}|^r<\infty$, via Theorems \ref
{prop3} and \ref{thm-next},
implies that $R_{n,\eps}<\infty$ for all $(n,\eps)\in\NN_{\geq
n_0}\times\RR^+$.
The proof of Theorem \ref{thm1} is hence completed by applying Lemmas
\ref{lem3} and \ref{lem4}.
\end{pf}
%
%%%%%%%%%%%%%%%%%%%%%%%%%5
%
\begin{lemma}\label{lem3}
Assume that there exists $n_0\geq2$ such that $R_{n,\eps}<\infty$
for all $(n,\eps)\in\NN_{\geq n_0}\times\RR^+$.
There then also exists $\delta>0$ such that $\lim_{n\to\infty
}R_{n,\delta}=0$.
\end{lemma}

%%%%%%%%%%%%%%%%%%%%%%%%%5
%
\begin{lemma}\label{lem4}
Statement \eqref{6'} is a consequence of Corollary \ref{cor_G-G-M}.
\end{lemma}
%
%%%%%%%%%%%%%%%%%%%%%%%%%%%%%%%%%%%%%%%%%%%%%%%%%%%%%%%%%%%%%%%%%%%
%
\begin{pf*}{Proof of Lemma \ref{lem3}}
%%%%%%%%%%%%%%%%%%%%%%%%%%%%%%%%%%%%%%%%%%%%%%%%%%%%%%%%%%%%%%%%%%%%%%%%%%%%%%%%%%%%%%%%%%%%%%%%%%%%%%%%%%%%%%%%%%%%%%%%%%%%%%%%%%%
We arrive at the conclusion from Lebesgue's dominated convergence
theorem, \cite{Cohn}, Theorem 2.4.4, page 72, by establishing that the
integrand
%
%e24 ###
\begin{equation}\label{6'd}
n^rh^{-(r+1)} \bigl( \bigl(\p(|X-x|< h |x| )
\bigr)^{n-1}-p_x^{n-1} \bigr)
\end{equation}
for some choice of $\delta$ and all $h\leq\delta$, is pointwise
decreasing in $n$ for sufficiently large $n$ and pointwise converging
to 0 as $n$ tends to infinity.
To this end, define
$\pi_x=\p(|X-x|< h |x| )$,
$g_x(y)=y^r (\pi_x^{y}-p_x^{y} )$, $ \lambda_1=-\log\pi_x$,
$\lambda_2=-\log p_x$.
To see that pointwise convergence to 0 holds, note that for some
$\delta$ and some $\eta>0$,
%
%e25 ###
\begin{equation}\label{6d}
\pi_x<1-\eta\qquad\mbox{for all }x \mbox{ and all } h<\delta.
\end{equation}
Condition (\ref{6d}) indeed prevails, except in the case where $F$ is
degenerate with total mass at a single point.
Given $\delta$ sufficiently small, $\pi_x^{n-1}-p_x^{n-1}$ therefore
decays exponentially in $n$, which yields
pointwise convergence to 0 of (\ref{6'd}).
The decreasing behaviour is equivalent to the existence of $y_0\geq0$
such that
%
%e26 ###
\begin{equation}
g_x(y_1) \geq g_x(y_2) \qquad\mbox{for all } y_1, y_2 \mbox{ such that }
y_0\leq y_1\leq y_2.\label{6''b}
\end{equation}
To verify \eqref{6''b}, note that
%
%e27 ###
\begin{equation}
g_x'(y) = -y^r (\lambda_1\mathrm{e}^{-\lambda_1 y}-\lambda_2\mathrm{e}^{-\lambda
_2 y} )
+ry^{r-1} (\mathrm{e}^{-\lambda_1 y}-\mathrm{e}^{-\lambda_2 y} ) = f_y(\lambda
_2)- f_y(\lambda_1) \label{6'b}
\end{equation}
with $f_y(\lambda):=\mathrm{e}^{-\lambda y} (\lambda y^r-ry^{r-1} )$
and furthermore that
%
%e28 ###
\begin{equation}
f_y'(\lambda) = \mathrm{e}^{-\lambda y}
(y^r-\lambda y^{r+1}+ry^r ) = \mathrm{e}^{-\lambda y}
\bigl((r+1)y^r-\lambda y^{r+1} \bigr). \label{6c}
\end{equation}
We verify (\ref{6''b}) using the fact that $f_y'(\lambda)<0$ for
$\lambda_1\leq\lambda\leq\lambda_2 $, which, by (\ref{6c}), is
satisfied for $y>y_0$, provided $ \lambda_1>\eta$ for some $\eta>0$.
The latter condition is equivalent to (\ref{6d}).
\end{pf*}

%%%%%%%%%%%%%%%%%%%%%%%%%%%%%%%%%%%%%%%%%%%%%%%%%%%%%%%%%%%%%%%%%%%

%%%%%%%%%%%%%%%%%%%%%%%%%%%%%%%%%%%%%%%%%%%%%%%%%%%%%%%%%%%%%%%%%%%
%
\begin{pf*}{Proof of Lemma \ref{lem4}}
It follows from Corollary \ref{cor_G-G-M} with $U_n=S^2_n/V^2_n$ that
%
%e29 ###
\begin{equation} \label{t1a}
\lim_{n\to\infty} \frac{r}{2}\int_0^n z^{r/2-1}\p(U_n>z
)\,\mathrm{d}z = \E|T|^r \qquad\mbox{for all $r>0$.}
\end{equation}
Define
$E_n= \{X_1=X_2=\cdots=X_n\neq0 \}$
so that
$\p(U_n>z )=\p(U_n^*>z )+\p(E_n)$ for $ 0<z<n$.
The desired conclusion is hence established by showing that for all $r>0$,
%
%e32 ###
%e31 ###
%e30 ###
\begin{eqnarray}\label{t1'b}
\lim_{n\to\infty} n^{r/2+1}\int_0^{n-\delta} z^{r/2-1}\p
(E_n ) (n-z)^{-(r/2+1)}\,\mathrm{d}z &=& 0,\\\label{t1b}
\lim_{n\to\infty} \int_0^{n-\delta} z^{r/2-1}\p(U_n>z )
\bigl(n^{r/2+1}(n-z)^{-(r/2+1)}-1 \bigr)\,\mathrm{d}z &=& 0,\\\label{t1c}
\lim_{n\to\infty} \int_{n-\delta}^n z^{r/2-1}\p(U_n>z
)\,\mathrm{d}z &=& 0.
\end{eqnarray}
%
%%%%%%%%%%%%%%%%%%%%%%%%%%%%%%%%%%%%%%%%%%%%%%%%%%%%%%%%%%%%%%%%%%%%%%%%%%%%%%
Starting with (\ref{t1'b}), let $\{a_k\}_{k\geq1}$ be a denumeration
of all non-zero points attributed mass by $F$ and define $p_k=\p
(X=a_k)$, $p=\sup_{k\geq1}p_k$.
It follows that $p<1$ since $X$ is not constant. Moreover,
\[
\p(E_n) = \sum_{k\geq1}p_k^n \leq p^{n-1}\sum_{k\geq1}p_k \leq p^{n-1}.
\]
This shows that $\p(E_n)$ decays exponentially in $n$. However, the quantities
\[
n(n-1)^{r/2}\int_0^{n-\delta} z^{(r-2)/2} (n-z)^{-(r+2)/2}\,\mathrm{d}z
\]
are all finite and grow with polynomial rate as $n$ grows. Conclusion
\eqref{t1'b} follows.
%%%%%%%%%%%%%%%%%%%%%%%%%%%%%%%%%%%%%%%%%%%%%%%%%%%%%%%%%%%%%%%%%%%%%%%%%%%%%%
Statement (\ref{t1c}) may be deduced from (\ref{t1a}) in the
following way:
\[
\int_{n-\delta}^n z^{r/2-1}\p(U_n>z )\,\mathrm{d}z
\leq(n-\delta)^{-1}\int_{n-\delta}^n z^{r/2}\p(U_n>z
)\,\mathrm{d}z \leq(n-\delta)^{-1}C_{r+2},
\]
where the constant $C_{r+2}$ stems from the identity in (\ref{t1a})
with $r$ replaced by $r+2$.
%%%%%%%%%%%%%%%%%%%%%%%%%%%%%%%%%%%%%%%%%%%%%%%%%%%%%%%%%%%%%%%%%%%%%%%%%%%%%%
It remains to prove (\ref{t1b}), which we split into
%
%e34 ###
%e33 ###
\begin{eqnarray} \label
{t1'd}
\lim_{n\to\infty} \int_0^{1} z^{(r/2-1)}\p(U_n>z )
\bigl(n(n-1)^{r/2}(n-z)^{-(r/2+1)}-1 \bigr)\,\mathrm{d}z&=&0,\\\label{t1d}
\lim_{n\to\infty} \int_1^{n-\delta} z^{r/2-1}\p(U_n>z )
\bigl(n(n-1)^{r/2}(n-z)^{-(r/2+1)}-1 \bigr)\,\mathrm{d}z&=&0.
\end{eqnarray}
Statement (\ref{t1'd}) follows from Lebesgue's dominated convergence
theorem, \cite{Cohn}, Theorem 2.4.4, page 72. To verify (\ref{t1d}),
we introduce the notation
\begin{eqnarray*}
f_n(z)&=&z^{r/2-1}\p(U_n>z )
\bigl(n(n-1)^{r/2}(n-z)^{-(r/2+1)}-1 \bigr)I_{D_n}, \\
D_n&=& \{z\dvtx  1\leq z\leq(n-\delta) \},\qquad
g_n(z)=z^{r}\p(U_n>z )I_{D_n},\qquad
g(z)=z^{r}\p(T^2>z )I_{D_n}.
\end{eqnarray*}
The desired conclusion (\ref{t1d}) is now written as (\ref{t1d''}),
while (\ref{t1d'}) follows from the assumptions, (\ref{t1a}) and the
elementary inequalities (\ref{t1d'''}):
%
%e37 ###
%e36 ###
%e35 ###
\begin{eqnarray}\label{t1d'''}
(n-1)/\bigl(z(n-z)\bigr) &\leq&(n-1)/\bigl(\delta(n-\delta)\bigr) \leq C \qquad\mbox{when }
z\in D_n,\\ \label{t1d''}
\lim_{n\to\infty}\int f_n&=&0,
\end{eqnarray}\vspace*{-12pt}
\begin{equation}\label{t1d'}
\int g_n\to\int g,\qquad g_n \to g,\qquad f_n\to0,\qquad |f_n|\leq C_1g_n.
\end{equation}
By a technique called \textit{Pratt's lemma}, Fatou's lemma, \cite
{Cohn}, Theorem 2.4.3, page~72, and (\ref{t1d'}) then give
%
%e39 ###
%e38 ###
\begin{eqnarray}\label{t1f}
C_1\int g&=& \int\liminf_{n}(C_1g_n-f_n)
\leq \liminf_{n} \int(C_1g_n-f_n)=C_1\int g-\limsup_{n} \int f_n,\qquad
\\
C_1\int g&=& \int\liminf_{n}(C_1g_n+f_n)\label{t1g}
\leq \liminf_{n} \int(C_1g_n+f_n)=C_1\int g+\liminf_{n} \int f_n.\qquad
\end{eqnarray}
Statement (\ref{t1d''}) follows from (\ref{t1f}) and (\ref{t1g}).
\end{pf*}

%%%%%%%%%%%%%%%%%%%%%%%%%%%%%%%%%%%%%%%%
\section*{Acknowledgements}
I would like to thank my Ph.D. supervisor Allan Gut for guidance,
encouragement and persistent reading of drafts.
I also wish to express my gratitude to Professor Lennart Bondesson for offering
valuable comments and criticism regarding \cite{jonsson} in connection
with the defense of my Licentiate thesis.

\printhistory

\end{document}